\documentclass[12pt]{article}
\usepackage{amsmath}
\usepackage{amssymb}
\usepackage{amsfonts}
\usepackage{amstext}
\usepackage{amsthm}

\DeclareMathOperator{\Iso}{Iso}
\DeclareMathOperator{\Aut}{Aut}
\DeclareMathOperator{\Gal}{Gal}
\DeclareMathOperator{\id}{id}
\DeclareMathOperator{\Tr}{Tr}
\DeclareMathOperator{\GL}{GL}
\DeclareMathOperator{\Ort}{O}
\DeclareMathOperator{\SL}{SL}
\DeclareMathOperator{\detp}{det}

\newcommand{\R}{{\mathbb R}}
\newcommand{\F}{{\mathbb F}}
\newcommand{\Q}{{\mathbb Q}}
\newcommand{\C}{{\mathbb C}}
\newcommand{\Z}{{\mathbb Z}}
\newcommand{\cL}{{\mathcal L}}

\newcommand{\q}{\quad}
\newcommand{\wt}{\widetilde}

\theoremstyle{plain}
\newtheorem{thm}{Theorem}[section]

\theoremstyle{remark}

\title{On the involutions fixing the class of a lattice}
\author{H.-G. Quebbemann, E.~M.~Rains}
\date{}

\begin{document}

\maketitle

\begin{abstract}
With any integral lattice $\Lambda$ in $n$-dimensional euclidean
space we associate an elementary abelian 2-group $I(\Lambda)$
whose elements represent parts of the dual
lattice that are similar to $\Lambda$.
There are corresponding involutions on modular forms for which the
theta series of $\Lambda$ is an eigenform;
previous work has focused on this connection.
In the present paper $I(\Lambda)$ is considered as a quotient of some
finite 2-subgroup of $\Ort_n(\R)$.
We establish upper bounds, depending only on $n$, for the order of
$I(\Lambda)$, and we study
the occurrence of similarities of specific types.
\end{abstract}

\section*{Introduction}

The arithmetic of lattices in euclidean spaces contains a counterpart to
the action of Atkin-Lehner involutions on elliptic modular forms.
An involution here takes the isometry class of a lattice to that of a rescaled
partial dual.
There are remarkable cases of (non-unimodular) lattices that are invariant under
this action; such a lattice is called strongly modular.
Some prominent examples have arisen as sections of the Leech lattice or in the
study of finite rational matrix groups.
See \cite{N}, \cite{Q}, \cite{R-S}, \cite{S-S}.

A basic question that has remained open is
when in general strong modularity can occur.
For instance, given $n=4r$ and $\ell =p_1\cdot \ldots \cdot p_s$
with distinct primes $p_i$ ($s>0$), there exist even
lattices of dimension $n$ and level $\ell$ whose genera have the necessary
property of being invariant under all $2^s$ involutions. But as we shall
see, a class can have this property only if $s\le 2 v_2(n)$, where $v_2$
denotes the 2-adic exponent valuation.
More generally, we find that at most $4^{[n/2]}$ involutions fix the class of a
given $n$-dimensional integral lattice. In cases of squarefree level here $[n/2]$
can be replaced by $v_2(n)$. These bounds are attained for any $n$ by certain
lattices built up from the case $n=2$ (which was studied in detail in \cite{R1}).

When one considers a single involution, again the question arises what additional
conditions must an invariant genus satisfy in order to contain an invariant class.
Some restrictions have been found by studying congruence properties mod 2 of the
relevant Atkin-Lehner eigenforms (\cite{R-S}, Appendix, and \cite{R2}). A more
elementary discussion related to this question is given in the final section of
this paper.

\section{Involutions and groups}

We shall study the involutory invariance of an integral lattice by applying
representation theory to its ``isomorphism group''. This object, which is defined below,
has occurred in \cite{N} as a subset of $\GL_n(\R)$, and as a group it had a short
appearance in \cite[Theorem 6]{R-S}.

We fix positive integers $\ell$ and $n$. Let $\cL$ denote the set of lattices $\Lambda$
in euclidean $\R^n$ such that $\Lambda\subseteq \Lambda^*$ and $\ell\Lambda^*\subseteq
\Lambda$, where $\Lambda^*$ is the dual lattice.
Then the orthogonal group $\Ort_n(\R)$ acts on $\cL$. Furthermore we have the
elementary abelian 2-group $W(\ell)$ of positive divisors $m\|\ell$ (i.e., $m|\ell$ and
$gcd(m,\ell/m) =1$), with multiplication $m*m' = mm'/gcd(m,m')^2$. Defining (as in
\cite[Section~1]{Q})
\[
   \Lambda_m = \sqrt{m}(\Lambda^* \cap m^{-1} \Lambda)
\]
we also have a group action of $W(\ell)$ on $\cL$, which clearly commutes with the
action of $\Ort_n(\R)$. In particular, either group acts on the orbits of the other one.

Now fix $\Lambda$, and let $\ell$ be its level in the sense that $\ell$ is minimal for
$\ell \Lambda^* \subseteq \Lambda$ (so $\ell$ divides $\det \Lambda = \# \Lambda^*/\Lambda$) .
 Then the following groups are canonically associated with $\Lambda$. First, $\Ort_n(\R)$
contains the usual stabilizer $\Aut(\Lambda)$ and the stabilizer of the $W(\ell)$-orbit
\[
   \Iso(\Lambda) = \{\sigma\in \Ort_n(\R) \mid \Lambda = \sigma\Lambda_m
                     \q\text{for some}\q m\in W(\ell)\}.
\]
Next, in $W(\ell)$ there are again the stabilizer of $\Lambda$ and the stabilizer of
its $\Ort_n(\R)$-orbit, denoted by $A(\Lambda)$ and $I(\Lambda)$, respectively.
Note that if $m \in I(\Lambda)$, then $\det \Lambda = \det \Lambda_m$, which implies
$m^{n/2} \| \det \Lambda$. Actually
$A(\Lambda)$ is trivial. In fact, if
$\Lambda = \Lambda_m$, then $m$ must be a square, say
$m=k^2$, and $(1/\sqrt{k})\Lambda$ is an integral lattice of level
$\ell/k$ and determinant prime to $m$; so $m=1$. As an immediate
consequence we obtain a group isomorphism
\[
   \Iso(\Lambda)/\Aut(\Lambda) \cong I(\Lambda)
\]
induced by $\sigma\longmapsto m: \Lambda = \sigma\Lambda_m$. Note also that the
squares in $I(\Lambda)$ form a subgroup $B(\Lambda)$ which corresponds to the
quotient of $\Iso(\Lambda)\cap \GL(\Lambda\otimes\Q)$ by $\Aut(\Lambda)$.
Of course, $B(\Lambda)= I(\Lambda)$ for odd $n$.

If $I(\Lambda) = W(\ell)$, then $\Lambda$ is said to be strongly modular.
Two ``trivial'' constructions of such lattices will be used. Namely, if $\Lambda$
and $\Lambda'$ are strongly modular with coprime levels $\ell,\ell'$, then
so are, with level $\ell\ell'$,
\begin{itemize}
\item[--]
the orthogonal sum $\sqrt{k'}\Lambda\oplus \sqrt{k}\Lambda'$ in the case
of squares $\ell = k^2$, $\ell'= k^{'2}$
\item[--]
the tensor product $\Lambda\otimes \Lambda'$.
\end{itemize}

\section{A general bound}

To measure $I(\Lambda)$ we first replace the linear group $\Iso(\Lambda)$ by a
2-Sylow subgroup $G$, so that $I(\Lambda) \cong G/P$ for $P = G\cap \Aut(\Lambda)$.
If $n$ is odd, an even larger elementary abelian quotient of $G$ arises from $P \cap
\SL_n(\R)$ because $P$ contains $-\id$. Put
\[
   d(G) = d(G/Q)
\]
for any finite $p$-group $G$, where $Q$ denotes the minimal normal subgroup such that
$G/Q$ is elementary abelian, in which case $d$ means dimension over $\F_p$. We would like
to use that $d(G) \le n$ if $G$ is a linear group of degree $n$. Actually this inequality
is true in general for odd $p$, while for 2-subgroups of $\GL_n(\C)$ one has $[3n/2]$ as
the best possible general upper bound (\cite{I}, \cite{W}).

\begin{thm} \label{thm1}
If $G$ is a finite $2$-subgroup of $\GL_n(\R)$, then $d(G) \le n$.
\end{thm}

This result can be extracted from \cite[Lemma~4]{W}, in the case where the
natural representation of $G$ on $\R^n$ is monomial, but the method equally
applies when this representation is induced from a 2-dimensional dihedral
representation of a subgroup. (The inductive proof in any case requires to
make a stronger statement, implying that if $H$ is a normal subgroup of
$G$, then $d(H) \le n$ and $d(G/H) \le n$ where at least one inequality is
strict. Of course, over $\C$ this already fails for $n=1$.) The remaining
cases of reducible or not absolutely irreducible groups are easily dealt
with, using induction on $n$ or an upper bound over $\C$ in dimension
$\nu = n/2$ (the crude bound $2\nu-1$ suffices).

\begin{thm} \label{thm2}
For any integral lattice $\Lambda \subset \R^n$,
\[
   d(I(\Lambda)) \le 2[n/2].
\]
In particular, if $\Lambda$ is strongly modular, at most
$2[n/2]$ distinct primes divide the level $\ell$. Lattices that attain this
bound exist in  any dimension $n$.
\end{thm}

\begin{proof}
The bound is an immediate consequence of the preceding theorem and the
remarks made at the beginning of this section. It is attained by
strongly modular lattices that arise by the first construction mentioned at
the end of Section~1. All we must know is that, for $n=2$ and primes
$p\not= q$, there exists such a lattice with level $p^2q^2$; for this see
\cite[Corollary~10]{R1}. For general even $n$ we then take scaled
orthogonal sums of such lattices, and for odd $n$ we add a 1-dimensional
summand.
\end{proof}

It should be noted at this point that also many pairs of primes $p\not= q$
admit a strongly modular 2-dimensional lattice with level $pq$; the
condition is that $p$ and $q$ must be quadratic residues of each other
(\cite{Q}, \cite{R1}). Taking now tensor products of such lattices, the
bound above is still attained for $n=4$, but otherwise this construction
suggests a different bound which will be obtained in  the next section.

It may be expected, however, that orthogonally indecomposable lattices
attaining the bound of Theorem~\ref{thm2} exist in any dimension. To obtain
an example for $n=3$, consider a solution of $a^2 + b^2 = c^2$ where
$a,b,c$ are coprime positive integers with $a+c$ divisible by 4. Take the
orthogonal sum $(a) \oplus \left(\begin{smallmatrix} c & b\\b &
c\end{smallmatrix}\right)$, and construct its Kneser
2-neighbour with respect to the sum of the first two base vectors. For this
second lattice $\Lambda$, which like the first has level $\ell = a^2$, we
find the Gram matrix
\[
   C =
  \begin{pmatrix}
   \frac{1}{4}(a+c) & \frac{b}{2} & \frac{1}{2}(c-a)\\
   \frac{b}{2} & c & b\\
   \frac{1}{2}(c-a) & b & a+c
   \end{pmatrix}
   = S(\ell C^{-1}) S^t, \qquad S = \begin{pmatrix}
                                    0 &  0 & 1\\
                                    0 & -1 & 0\\
                                    1 &  0 & 0
                                    \end{pmatrix}\, .
\]
This shows that again $\{1,\ell\} \subseteq I(\Lambda)$. For $(a,b,c) = (15,8,17)$, computation of $\Lambda_9$ moreover
shows that $\Lambda$ is strongly modular; it is clearly indecomposable. ($C$ is reduced in the sense of
Seeber-Eisenstein \cite{O} whenever $a = 4t^2-1$ and $b=4t$, $t>1$. For $t>2$, however, $I(\Lambda)$ does not contain
$m=(2t-1)^2 \in W(\ell)$. In fact, $\Lambda$ has minimal norm $2t^2$, while $\Lambda_m$ contains the vector of norm 8
given by $(1/\sqrt{m})(2,-2,1)$ with respect to $C$.)

\section{The squarefree case}

If we only consider lattices $\Lambda$ with squarefree level or, in the general case,
factor out $B(\Lambda) = I(\Lambda) \cap \{k^2\mid k=1,2,\ldots\}$, then the result of
Section~2 can be improved considerably. We now make use of the fact that any element of
$\Iso(\Lambda)$ up to some scalar factor $1/\sqrt{m}$ is defined over $\Z$, and so
$\Iso(\Lambda)$ is conjugate to a subgroup of
\[
   \GL_n(\Q)^+ = \langle \GL_n(\Q), \; \{\sqrt{p}\cdot\id \mid p
                \text{ prime}\}\rangle \, .
\]

\begin{thm} \label{thm3}
Let $G$ be a finite $2$-subgroup of $\GL_n(\Q)^+$, and define $N = G\cap \GL_n(\Q)$. Then
$d(G/N) \le 2 v_2(n)$.
\end{thm}

\begin{proof}
We can write $G = G_1 \cup G_2 \cup \ldots$, where $G_m = G\cap \sqrt{m} \GL_n(\Q)$ is a
coset of $N=G_1$ or is empty. We may assume that $N$ contains $-\id$ (or adjoin it). Let
$\varrho$ denote the natural representation of $G$ on $V = K^n$, where $K$ is the field
obtained from $\Q$ by adjoining the relevant square roots $\sqrt{m}$ (i.e., those for
which $G_m$ is nonempty) together with $\zeta = e^{2\pi i/r}$, where $r$ is the exponent
of $G$. Let $\chi$ be the character of $\varrho$. Then the values of $\chi$ lie in the
2-power cyclotomic field $\Q(\zeta)$, while on the other hand $\chi(\sigma) \in
\Q\sqrt{m}$ for $\sigma \in G_m$. This shows that $\chi$ vanishes outside $N\cup G_2$.
Set $M = N$ if $\chi$ also vanishes on $G_2$, and $M = N\cup G_2$ otherwise. Denoting the
restriction of $\chi$ to the subgroup $M$ by $\chi_M$, we obtain
\[
   \langle \chi,\chi\rangle
    = \frac{1}{\# G} \sum_{\sigma\in G} |\chi(\sigma)|^2
    = \frac{1}{\# G} \sum_{\sigma\in M} |\chi(\sigma)|^2
    = \frac{1}{(G:M)} \langle \chi_M, \chi_M\rangle\, .
\]
It follows that $(G:M) \le n^2$, since $\langle \chi,\chi\rangle \ge 1$ and
$\langle \chi_M,\chi_M\rangle \le n^2$. However, the theorem moreover
asserts that $(G:N)$ divides $n^2$. To prove this we first observe that the
Galois group $\Gamma = \Gal(K/ \Q)$ acts on the $G$-invariant subspaces of
$V$. If such a subspace $U\not= 0$ is defined over $\Q$ (i.e., invariant
also under $\Gamma$), then the kernel of the restriction of $\varrho$ to
$U$ clearly is contained in $N$, and in the case $U\not= V$ it suffices to
apply the theorem on $U$ and $V/U$ (using induction on $n$). Therefore we
now assume that no such proper subspace exists. This implies that the
isotypical summands of $\varrho_M$ are permuted transitively by the
(semidirect) product $G\cdot \Gamma$. So we have
\[
   \chi_M = a(\psi_1 + \cdots + \psi_b)
\]
where $\psi_1,\ldots,\psi_b$ are pairwise inequivalent irreducible
characters of $M$ and all are of the same degree $c$. Then $n= abc$ and
\[
   n^2 = bc^2 \langle \chi_M, \chi_M\rangle = bc^2 \langle \chi,
         \chi\rangle (G:M)\, .
\]
Finally suppose that $(M:N) = 2$. Then it suffices to show that $bc^2$ is
even. Since $M$ is a 2-group, $c$ is a power of 2, and we only have to show
that $b$ is even if $c=1$. Let $\psi$ be one of the linear characters
$\psi_i$. If $\psi^2 = 1$, set $\wt{\psi} = \lambda\psi$ where $\lambda$
is the linear character of $M$ with kernel $N$, and if $\psi^2 \not= 1$ set
$\wt{\psi} = \psi^{-1}$. It follows that $\wt{\psi} \not= \psi$ and
$\wt{\psi}$ also appears in $\chi_M$ (to see this in the first case, use a
Galois automorphism mapping $\sqrt{2}$ to $-\sqrt{2}$). Therefore the
linear characters $\psi_i$ appear in pairs.
\end{proof}

\begin{thm} \label{thm4}
For any integral lattice $\Lambda \subset \R^n$,
\[
   d(I(\Lambda)/B(\Lambda)) \le 2v_2(n)\, .
\]
In particular, if $\Lambda$ is strongly modular with squarefree level $\ell$, then $\ell$
is the product of at most $2v_2(n)$ primes. Lattices that attain this bound exist in any
dimension $n$.
\end{thm}

\begin{proof}
Recall that $I(\Lambda)/B(\Lambda) \cong \Iso(\Lambda) / \Iso(\Lambda)\cap
\GL(\Lambda\otimes \Q)$, where the last two groups again may be replaced by Sylow
2-subgroups. So the bound immediately follows from the preceding theorem. For $n=2^e f$,
$f$ odd, this bound is attained by the $f$-fold orthogonal sum of the tensor product of
$e$ strongly modular 2-dimensional lattices with coprime levels $p_i q_i$ ($p_i \not=
q_i$ primes). As was noted already, enough such lattices exist.
\end{proof}

If more is known on the action of a Sylow 2-subgroup of $\Aut(\Lambda)$,
then of course the proof of Theorem~\ref{thm3} gives more information on
$I(\Lambda)$. Suppose, for example, that this action is multiplicity-free.
Then in that proof the equation $n = abc$ holds with $a=1$ and $b = \langle
\chi_M,\chi_M\rangle = (G:M) \langle \chi,\chi\rangle$, so for $M=N$ (in
particular, for odd $\ell$) we obtain
\[
   d(I(\Lambda) / B(\Lambda)) \le v_2(n)\, .
\]
Given $\ell = p_1\cdot \ldots \cdot p_s$ with $s$ distinct primes $p_i$, this bound is attained by the tensor product
of 2-dimensional lattices having the Gram matrix $\left(\begin{smallmatrix} 1 & 0\\0 & p_i\end{smallmatrix}\right)$ or
$\left(\begin{smallmatrix} 2 & 1\\1 & q_i\end{smallmatrix}\right)$, $q_i = (p_i+1)/2$. The condition of multiplicity
one is satisfied because each of these lattices admits a reflection and the corresponding automorphisms of $\Lambda$
generate an elementary abelian 2-group acting by its regular representation.

Moreover, it is this last bound, $s \le v_2(n)$ for a strongly modular
lattice, which some of the ``extremal'' lattices discovered during the
past few years attain.

\section{Single involutions}

In this section we fix a positive integer $m$ which is not of the form
$k^2$ or $2k^2$ for an integer $k$.
If $m$ belongs to $I(\Lambda)$
for an integral lattice $\Lambda$, then
as before $\Lambda$ is the image of $\Lambda^* \cap m^{-1}\Lambda$
under a similarity
$\tau=\sqrt{m}\sigma$ where $\sigma \in \Ort_n(\R)$ is of 2-power order;
that is, we have $f(\tau) =0$ with $f= X^{2r}-m^r$ for some $r=2^k$, $k \geq 0$.
Let such a polynomial $f$ be fixed, and define the algebra
\[
E = \Q[X]/(f) = \Q[\alpha], \quad \alpha = X+(f).
\]
In $E$ we consider the root of unity $\zeta = \alpha^2/m$ and the order $R =
\Z[\alpha,\zeta]$ (over which our lattices will be defined). Now $f$
factors into the rationally irreducible polynomials
\[
f_0 = X^2-m, \quad f_j = X^{2^j}+m^{2^{j-1}}  (1 \leq j \leq k).
\]
So we have $E \cong E_0 \times ... \times E_k$ where $E_j
= \Q[X]/(f_j) =\Q(\alpha_j)$ is a quadratic extension of the $2^j$-th
cyclotomic field $\Q(\zeta_j)$, with $\alpha_j^2 =m\zeta_j$. Defining
$\overline{\alpha} = m \alpha^{-1}$ we obtain an involutory automorphism
"bar" on $E$, and correspondingly on each $E_j$. For $j\geq 1$ this involution
of $E_j$ is not the identity, so then $E_j$ also is
a (totally complex) quadratic extension of the totally real field
$\Q(\beta_j)$, $\beta_j=\alpha_j +\overline{\alpha_j}$, with $\beta_j^2
=m(1+\zeta_j)(1+\zeta_j^{-1})$. Note that $R = \overline{R}$.

By an inner product space over $E$, we mean a
finitely generated $E$-module $V$ carrying an $E$-valued, totally
positive definite hermitian form $h$ for the involution defined above.
Clearly such a hermitian module splits as an orthogonal sum of hermitian
spaces over the fields $E_j$ (for $j=0$ "hermitian" here means "symmetric bilinear").
By definition, $h$ is totally positive definite if each of these summands is.
In this case $V$ also  carries the rational inner product
\[
   u\cdot v = (2r)^{-1} \Tr_{E/\Q} h(u,v)
\]
where $\Tr_{E/\Q}$ denotes the sum of the field traces.
Given an $R$-lattice $\Lambda \subset V$, we denote its $R$-dual
with respect to $h$ by $\Lambda^h$ and its $\Z$-dual with respect to
the dot product by $\Lambda^*$.

\begin{thm} \label{thm5}
If $(V,h)$ is an inner product space over $E$ and $\Lambda$ is an
$R$-lattice on $V$, then the duals are related by $\Lambda^h=\tau(\Lambda^*)$,
where $\tau(v)=\alpha v$.
In particular, $\Lambda$ is contained in $\Lambda^h$ with index prime to $m$
precisely when
\[
\Lambda \subseteq \Lambda^*, \quad m \in I(\Lambda),
\quad \Lambda=\tau(\Lambda^*\cap m^{-1}\Lambda).
\]
Conversely, if these conditions
are satisfied by a lattice $\Lambda \subset \R^n$ and a similarity $\tau$ with
$f(\tau)=0$, then $\Lambda$ and $\tau$ arise from $R$ in the way above.
\end{thm}

\begin{proof}
For the $\Z$-basis
$1,\zeta,...,\zeta^{r-1},\alpha,\alpha\zeta,...,\alpha\zeta^{r-1}$ of $R$
it is easily seen that all elements have trace 0, except for $\Tr_{E/\Q}(1)=2r$.
Therefore, given $\varepsilon \in E$, the condition
$(2r)^{-1}\Tr_{E/\Q}(\varepsilon R) \subseteq \Z$ is equivalent to
$\varepsilon \in \alpha^{-1}R$. This proves the first part.
As to the converse, the action of $\tau$ makes
$V = \Lambda \otimes \Q$ into an $E$-module, and since $\Lambda$ is
stable under both $\tau$ and $(1/m)\tau^2$, it is an $R$-lattice.
We have to show that the inner product arises from one over $E$ in the way
described. For fixed $u$ and $v$ in $V$ consider the linear map $t$ from
$E$ to $\Q$ defined by $t(\varepsilon)= (\varepsilon u)\cdot v$. There is
a unique $\varepsilon '$ in $E$ such that
\[
\Tr_{E/\Q} (\varepsilon \varepsilon') = 2r t(\varepsilon)
\]
for all $\varepsilon \in E$, and we define $h(u,v)$ to be this element
$\varepsilon'$. The verification of further details will be omitted.
(Compare \cite{M}; the discussion restricted there to isometries carries over
to similarities).
\end{proof}

In the following $m$ is supposed to be odd. The theorem may be applied to
find restrictions for the rational quadratic space
$V = \Lambda \otimes \Q$ when $I(\Lambda)$ contains $m$.
Recall that in this case $\det \Lambda = m^{n/2} d'$ where $m$ and $d'$ are
relatively prime.
The necessary condition that $V$ admits the similarity factor $m$ just
means that (i) $m$ is a square modulo each prime $p>2$ dividing $d'$ with odd
multiplicity and (ii) $d'$ is a square modulo each prime $p$ dividing $m$ with
odd multiplicity.
This can be seen, for instance, by \cite{S}, Chapter 5, Corollary 3.6.

Let $\Lambda$ be an even lattice. We use the rational invariants defined
at an arbitrary prime number $p$ by the Gauss sum
\[
\gamma_p(\Lambda) = (\detp _p \Lambda)^{-1/2} \sum \exp(\pi i v\cdot v),
\]
where the summation extends over the elements $v+\Lambda$ in the $p$-primary
component of $\Lambda^*/\Lambda$; $\det_p \Lambda$ denotes the order of this group.
The product $\prod \gamma_p(\Lambda)$ taken over all primes $p$ is $i^{n/2}$.
(See \cite{S}, Chapter 5, {\S} 8.)

Recall that if $m \in I(\Lambda)$, then
$V$ is defined over the direct product of
the fields $E_0 =\Q(\sqrt{m})$, $E_1 = \Q(\sqrt{-m})$, $E_2, \ldots, E_k,$
where $E_j$ for $j >1$ has degree $2^{j-1}$ over $\Q(i)$.
Assume for the moment that the contributions from $E_0$ and $E_1$ are zero.
Then the dimension $n$ is a multiple of 4, and
the rational quadratic form on $V$ arises from a hermitian form over
$\Q(i)$, so it is an orthogonal sum $\phi \perp \phi$ and has square determinant.
We compute that for $p^s \|m$,
\[
\gamma_p(\Lambda) =
\begin{cases}
1 & \text{if} \q p^s \equiv 1 \pmod{4}\\
(-1)^{n/4} & \text{if} \q p^s \equiv 3 \pmod{4},
\end{cases}
\]
which gives
\[
\prod_{p|m} \gamma_p(\Lambda) =
(-1)^{\frac{n}{4} \cdot \frac{m-1}{2}}.
\]
When $m \equiv 7 \pmod{8}$, the assumption used here is eliminated
by the following theorem (first proved in \cite{R2} by less elementary
arguments, a special case can be found in \cite{R-S}).

\begin{thm} \label{thm6}
Let $m$ be congruent to $7$ mod $8$, and let $m'$ be a
positive integer relatively prime to $m$ such that $-m$ is a square mod
$m'$.  Then for any even lattice $\Lambda$ of level dividing $mm'$ with
$m \in I(\Lambda)$,
$$
\prod_{p|m'} \gamma_p(\Lambda)
=
\prod_{p|m'}
\left( \frac{\sqrt{-m}}{\det_p\Lambda} \right),
$$
where in the factor corresponding to $p$, $\sqrt{-m}$ represents either
square root of $-m$ in the $p$-adic integers.
\end{thm}

\begin{proof}
First, note that for odd $p|m'$ such that $m$ is not a square modulo $p$,
$\det_p\Lambda$ must be a square (since otherwise $V$ does not admit the
similarity factor $m$).  In particular, since $-m$ is assumed to be a
square modulo $p$, we find that either $\det_p\Lambda$ is a square or
$-1$
is a square modulo $p$, and thus the given Legendre symbols are all
well-defined.

Now consider the $R$-lattice structure for the order $R$ as before.
Let $R_j$ be the image of $R$ in $E_j$ for $j=0, \dots, k$, and define the
order
$$R'= R_0 \times \ldots \times R_k.$$
Since the index of $R$ in $R'$ is a power of 2, there exists a
positive integer $t$ such that $R'(2^t\Lambda)$ is an even lattice satisfying the
hypotheses (possibly $m'$ has to be multiplied by a power of 2), with the same
Gauss sums and $p$-determinants modulo squares as $\Lambda$.  Since $R'$
is a product of orders preserved by the involution,
$R'(2^t\Lambda)$ splits as a product of lattices,
one over each order.  Since the desired identity is preserved under taking
products, we may reduce to the case that $\Lambda$ is a lattice over
one of the factors $R_j$ of $R'$.

Case $j>1$. As noted already, here $\Lambda$ has square determinant and
dimension divisible by 4. The computations preceding the theorem and the product
formula give
\[
\prod_{p|m'} \gamma_p(\Lambda) = (-1)^{n/4} \prod_{p|m} \gamma_p(\Lambda)=1
\]
as required.

Case  $j=1$, $\Q(\sqrt{-m})$.  The rational quadratic form on
$V$ now splits as $\phi\perp m\phi$; since $-m$ is a
square
in $\Z_p$, $V\otimes \Q_p$ is hyperbolic for all $p|m'$.  In
particular, $\det_p \Lambda$ is a square and $\gamma_p(\Lambda)=1$ for
all
$p|m'$.

Case $j=0$, $\Q(\sqrt{m})$.  Then $\Lambda$ contains a sublattice of
index prime to $m$ that splits as an orthogonal sum  of 2-dimensional lattices
preserved by $\Z[\sqrt{m}]$.  (Indeed, if $v\in\Lambda$ has
norm
prime to $m$, then the sum of $\Lambda \cap \Q(\sqrt{m})v$ and
$\Lambda \cap (\Q(\sqrt{m})v)^\perp$ has index prime to $m$; proceed by
induction).  It thus suffices to consider the case of a primitive
2-dimensional even lattice of level exactly $mm'$, in which case the
claim follows from Theorem 11 of \cite{R1}.
\end{proof}

For example, let $4|n$ and let $\ell > 1$ be squarefree.
Given signs $\gamma_p = \pm 1$ for the primes $p|\ell$ such that
$\prod \gamma_p = (-1)^{n/4}$, there exist even lattices $\Lambda$
of dimension $n$, level $\ell$, determinant $\ell^{n/2}$ and rational invariant
$\gamma_p$ at each $p$. In this situation all divisors of $\ell$
are similarity factors of $\Lambda \otimes \Q$. However, when $\ell = mm'$
where $m \equiv 7 \pmod{8}$ and $-m$ is a square modulo $m'$, we obtain
$$\prod_{p|m'} \gamma_p =1$$ as a necessary condition for having
$m \in I(\Lambda)$. For $m$ incongruent to 7 mod 8 there is no similar restriction
in general.

\bigskip
\begin{tabular}{l@{\hspace*{1cm}}l}
Fachbereich 6 Mathematik & Center for Communication Research\\
Universit\"at Oldenburg & Institute for Defense Analyses\\
26111 Oldenburg, Germany & Princeton, NJ 08540, U.S.A.\\
quebbemann@mathematik. & rains@idaccr.org\\
uni-oldenburg.de\\
\end{tabular}

\end{document}